\documentclass{amsart}
\newtheorem{theorem}{Theorem}[section]
\newtheorem{lemma}[theorem]{Lemma}

\newtheorem{proposition}[theorem]{Proposition}

\theoremstyle{definition}
\newtheorem{definition}[theorem]{Definition}

\theoremstyle{remark}

\newcommand{\C}{{\mathcal C}}

\newcommand{\cP}{{\mathcal P}}

\newcommand{\bC}{{\mathbb{C}}}
\newcommand{\bN}{{\mathbb{N}}}

\newcommand{\ov}[1]{\overline{#1}}
\numberwithin{equation}{section}

\newcommand{\Tr}{\mathrm{Tr}}

\newcommand{\jed}{{\mathbb{I}}}

\title[Positive maps between $M_2(\bC)$ and $M_n(\bC)$]{On
decomposability of positive maps between $M_2(\bC)$ and $M_n(\bC)$}
\author[W. A. Majewski]{W{\l}adys{\l}aw A. Majewski}
\address{Institute of Theoretical Physics and Astrophysics, Gda{\'n}sk University,
Wita Stwosza 57, 80-952 Gda{\'n}sk, Poland}
\email{fizwam@univ.gda.pl}
\author[M. Marciniak]{Marcin Marciniak}
\email{matmm@univ.gda.pl}
\date{\today}
\keywords{Positive maps, decomposable maps, face structure}
\subjclass[2000]{47B65, 47L07}
\thanks{Supported by KBN grant PB/1490/PO3/2003/25. The authors would like also thank the support of EU RTN HPRN-CT-2002-00729.}
\begin{document}
\begin{abstract}
A map $\varphi:M_m(\bC)\to M_n(\bC)$ is decomposable if it is of the
form $\varphi=\varphi_1+\varphi_2$ where $\varphi_1$ is a CP map
while $\varphi_2$ is a co-CP map. A partial characterization of decomposability for maps $\varphi: M_2(\bC) \to M_3(\bC)$ is given.
\end{abstract}

\maketitle
\section{Introduction}
Let $\varphi:M_m(\bC)\to M_n(\bC)$ be a linear map. We say that
$\varphi$ is positive if $\varphi(A)$ is a positive element in
$M_n(\bC)$ for every positive matrix from $M_m(\bC)$. If $k\in\bN$,
then $\varphi$ is said to be $k$-positive (respectively
$k$-copositive) whenever $[\varphi(A_{ij})]_{i,j=1}^k$ (respectively
$[\varphi(A_{ji})]_{i,j=1}^k$) is positive in $M_k(M_n(\bC))$ for
every positive element $[A_{ij}]_{i,j=1}^k$ of $M_k(M_m(\bC))$. If
$\varphi$ is $k$-positive (respectively $k$-copositive) for every
$k\in\bN$ then we say that $\varphi$ is completely positive
(respectively completely copositive). Finally, we say that the map
$\varphi$ is decomposable if it has the form
$\varphi=\varphi_1+\varphi_2$ where $\varphi_1$ is a completely
positive map while $\varphi_2$ is a completely copositive one.

By $\cP(m,n)$ we denote the set of all positive maps acting between
$M_m(\bC)$ and $M_n(\bC)$ and by $\cP_1(m,n)$ -- the subset of
$\cP(m,n)$ composed of all positive unital maps (i.e. such that
$\varphi(\jed)=\jed$). Recall that $\cP(m,n)$ has the structure of a
convex cone while $\cP_1(m,n)$ is its convex subset.

In the sequel we will use the notion of a face of a convex cone.
\begin{definition}\label{face}
Let $C$ be a convex cone. We say that a convex subcone $F\subset C$
is a face of $C$ if for every $c_1,c_2\in C$ the condition
$c_1+c_2\in F$ implies $c_1,c_2\in F$.

A face $F$ is said to be a maximal face if $F$ is a proper subcone
of $C$ and for every face $G$ such that $F\subseteq G$ we have $G=F$
or $G=C$.
\end{definition}

The following theorem of Kye gives a nice characterization of
maximal faces in $\cP(m,n)$.
\begin{theorem}[\cite{Kyecan}]\label{maxface}
A convex subset $F\subset\cP(m,n)$ is a maximal face of $\cP(m,n)$
if and only if there are vectors $\xi\in\bC^m$ and $\eta\in\bC^n$
such that $F=F_{\xi,\eta}$ where
\begin{equation}\label{F}
F_{\xi,\eta}=\{\varphi\in\cP(m,n):\,\varphi(P_\xi)\eta=0\}
\end{equation}
and $P_\xi$ denotes the one-dimensional orthogonal projection in
$M_m(\bC)$ onto the subspace generated by the vector $\xi$.
\end{theorem}

The aim of
this paper is to discuss the problem whether it is possible to find
concrete examples of the decomposition $\varphi=\varphi_1+\varphi_2$ onto
completely positive and completely copositive parts
where $\varphi\in\cP(m,n)$. 
It is well known (see \cite{St,Wor}) that every elements of
$\cP(2,2)$, $\cP(2,3)$ and $\cP(3,2)$ are decomposable. 
In \cite{MM2} we proved that if $\varphi$ is extremal element
of $\cP_1(2,2)$ then its decomposition is unique. Moreover, we provided
a full description of this decomposition. 
In the case $m>2$ or $n>2$ the problem of finding decomposition is still
unsolved. In this paper we consider the next step for solving this 
problem, namely for the case $m=2$ and $n=3$.
Our approach will be based on the method
of the so called Choi matrix. 

Recall (see \cite{Choi,MM2}) for
details) that if $\varphi:M_m\to M_n$ is a linear map and
$\{E_{ij}\}_{i,j=1}^m$ is a system of matrix units in $M_m(\bC)$,
then the matrix
\begin{equation}\label{Cm}
\mathbf{H}_\varphi=[\varphi(E_{ij})]_{i,j=1}^m\in
M_m(M_n(\bC)),\end{equation} is called the Choi matrix of $\varphi$
with respect to the system $\{E_{ij}\}$. Complete positivity of
$\varphi$ is equivalent to positivity of $\mathbf{H}_\varphi$ while
positivity of $\varphi$ is equivalent to block-positivity of
$\mathbf{H}_\varphi$ (see \cite{Choi,MM2}).
Recall (see Lemma 2.3 in \cite{MM2}) that in the case $m=n=2$ the general form of the Choi matrix
of a positive map $\varphi$ is the following
\begin{equation}
\label{Choi22}
\mathbf{H}_\varphi=\left[\begin{array}{cc|cc}a&c&0&y\\\overline{c}&b&\overline{z}&t\\\hline
0&z&0&0\\\overline{y}&\overline{t}&0&u\end{array}\right]
\end{equation}
where $a,b,u\geq 0$, $c,y,z,t\in\mathbb{C}$ and the following inequalities are satisfied:
\begin{enumerate}
\item $|c|^2\leq ab$,
\item $|t|^2\leq bu$,
\item $|y|+|z|\leq (au)^{1/2}$.
\end{enumerate}
It will turn out that in the case $n=3$ the Choi matrix has the form which 
similar to (\ref{Choi22}) but some of the coefficients have to be matrices.
The main result of our paper is the generalization of Lemma 2.3 from \cite{MM2}
in the language of some matrix inequalities. It is worth to pointing out that technical lemmas
leading to this generalization are formulated and proved for more general case, i.e. for
$\varphi:M_2(\bC)\to M_{n+1}(\bC)$ where $n\geq 2$.

\section{Main results}
In this section we will make one step further in the analysis of
positive maps and we will examine maps $\varphi$ in $\cP_1(2,n+1)$
where $n\geq 1$. Let $\{e_1,e_2\}$ and $\{f_1,f_2,\ldots,f_{n+1}\}$
denote the standard orthonormal bases of the spaces $\bC^2$ and
$\bC^n$ respectively, and let $\{E_{ij}\}_{i,j=1}^2$ and
$\{F_{kl}\}_{k,l=1}^{n+1}$ be systems of matrix matrix units in
$M_2(\bC)$ and $M_{n+1}(\bC)$ associated with these bases. We assume
that $\varphi\in F_{\xi,\eta}$ for some $\xi\in\bC^2$ and
$\eta\in\bC^{n+1}$. By taking the map $A\mapsto V^*\varphi(WAW^*)V$
for suitable $W\in U(2)$ and $V\in U(n+1)$ we can assume without
loss of generality that $\xi=e_2$ and $\eta=f_1$. Then the Choi
matrix of $\varphi$ has the form
\begin{equation}\label{Cm1}
\mathbf{H}=\left[\begin{array}{ccccc|ccccc}
a&c_1&c_2&\ldots&c_n&x&y_1&y_2&\ldots&y_n \\
\ov{c_1}&b_{11}&b_{12}&\ldots&b_{1n}&\ov{z_1}&t_{11}&t_{12}&\ldots&t_{1n}\\
\ov{c_2}&b_{21}&b_{22}&\ldots&b_{2n}&\ov{z_2}&t_{21}&t_{22}&\ldots&t_{2n}\\
\vdots &\vdots&\vdots& &\vdots&\vdots&\vdots&\vdots& &\vdots \\
\ov{c_n}&b_{n1}&b_{n2}&\ldots&b_{nn}&\ov{z_n}&t_{n1}&t_{n2}&\ldots&t_{nn}\\\hline
\ov{x}&z_1&z_2&\ldots&z_n&0&0&0&\ldots&0\\
\ov{y_1}&\ov{t_{11}}&\ov{t_{21}}&\ldots&\ov{t_{n1}}&0&u_{11}&u_{12}&\ldots&u_{1n}\\
\ov{y_2}&\ov{t_{12}}&\ov{t_{22}}&\ldots&\ov{t_{n2}}&0&u_{21}&u_{22}&\ldots&u_{2n}\\
\vdots&\vdots&\vdots& &\vdots&\vdots&\vdots&\vdots& &\vdots\\
\ov{y_n}&\ov{t_{1n}}&\ov{t_{2n}}&\ldots&\ov{t_{nn}}&0&u_{n1}&u_{n2}&\ldots&u_{nn}
\end{array}\right]
\end{equation}
We introduce the following notations:
$$C=\left[\begin{array}{cccc}c_1&c_2&\ldots&c_n\end{array}\right],\quad
Y=\left[\begin{array}{cccc}y_1&y_2&\ldots&y_n\end{array}\right],\quad
Z=\left[\begin{array}{cccc}z_1&z_2&\ldots&z_n\end{array}\right],$$
$$B=\left[\begin{array}{cccc}b_{11}&b_{12}&\ldots&b_{1n}\\b_{21}&b_{22}&\ldots&b_{2n}\\\vdots&\vdots&
&\vdots\\b_{n1}&b_{n2}&\ldots&b_{nn}\end{array}\right],\quad
T=\left[\begin{array}{cccc}t_{11}&t_{12}&\ldots&t_{1n}\\t_{21}&t_{22}&\ldots&t_{2n}\\\vdots&\vdots&
&\vdots\\t_{n1}&t_{n2}&\ldots&t_{nn}\end{array}\right],\quad
U=\left[\begin{array}{cccc}u_{11}&u_{12}&\ldots&u_{1n}\\u_{21}&u_{22}&\ldots&u_{2n}\\\vdots&\vdots&
&\vdots\\u_{n1}&u_{n2}&\ldots&u_{nn}\end{array}\right].$$ 
The matrix
(\ref{Cm1}) can be rewritten in the following form
\begin{equation}\label{Cm2}
\mathbf{H}=\left[\begin{array}{cc|cc} a&C&x&Y \\ C^*&B&Z^*&T
\\\hline \ov{x}&Z&0&0 \\ Y^*&T^*&0&U\end{array}\right].
\end{equation}
The symbol $0$ in the right-bottom block has three different
meanings. It denotes
$
0
$, $\left[\begin{array}{cccc}0&0&\ldots&0\end{array}\right]$ or
$\left[\begin{array}{c}0\\0\\\vdots\\0\end{array}\right]$
respectively. We have the following
\begin{proposition}\label{Cmpos}
Let $\varphi:M_2(\bC)\to M_{n+1}(\bC)$ be a positive map with the
Choi matrix of the form (\ref{Cm2}). Then the following relations
hold:
\begin{enumerate}
\item $a\geq 0$ and $B$, $U$ are positive matrices,
\item if $a=0$ then $C=0$, and if $a>0$ then $C^*C\leq aB$,
\item $x=0$,
\item the matrix
$\left[\begin{array}{c|c}B&T\\\hline T^*&U\end{array}\right]\in
M_2(M_n(\bC))$ is block-positive.
\end{enumerate}
\end{proposition}
\begin{proof}
It follows from positivity of $\varphi$ that blocks on main
diagonal, i.e.
$\varphi(E_{11})=\left[\begin{array}{cc}a&C\\C^*&B\end{array}\right]$
and
$\varphi(E_{22})=\left[\begin{array}{cc}0&0\\0&U\end{array}\right]$,
must be positive matrices. This immediately implies (1) and (2) (cf.
\cite{Wor}). From block-positivity of $\mathbf{H}$ we conclude that
the matrix
$$\left[\begin{array}{cc}a&x\\\ov{x}&0\end{array}\right]=
\left[\begin{array}{cc}\langle f_1,\varphi(E_{11})f_1\rangle &
\langle f_1,\varphi(E_{12})f_1\rangle \\ \langle
f_1,\varphi(E_{21})f_1\rangle & \langle
f_1,\varphi(E_{22})f_1\rangle\end{array}\right]$$ is a positive
element of $M_2(\bC)$. So $x=0$, and (3) is proved. The statement of
point (4) is an obvious consequence of block-positivity of
$\mathbf{H}$.
\end{proof}

For $X=\left[\begin{array}{cccc}x_1&x_2&\ldots&x_n\end{array}\right]\in M_{1,n}(\bC)$ we define $\Vert X\Vert=\left(\sum_{i=1}^n|x_i|^2\right)^{1/2}$.
By $|X|$ we denote the square ($n\times n$)-matrix $(X^*X)^{1/2}$.
Let us observe that for any $X\in M_{1,n}(\bC)$ we have
$$
|X|=\Vert X\Vert P_{\xi_X}
$$
where $\xi_X=\Vert X\Vert^{-1}X^*$ and $P_\xi$ denotes the orthogonal 
projection onto the one-dimensional subspace in $\bC^n$ generated by
a vector $\xi\in\bC^n$ (we identify elements of $M_{n,1}(\bC)$ with vectors from $\bC^n$).
We have the following
\begin{lemma}\label{lemmain}
A map $\varphi$ with the Choi matrix of the form
\begin{equation}\label{Cm3}
\mathbf{H}=\left[\begin{array}{cc|cc} a&C&0&Y \\ C^*&B&Z^*&T
\\\hline 0&Z&0&0 \\ Y^*&T^*&0&U\end{array}\right].
\end{equation}
is positive if and only if the inequality
\begin{equation}\label{nier1}
\left|\langle Y^*,\mathit{\Gamma}^\tau\rangle+
\ov{\langle Z^*,\mathit{\Gamma}^\tau\rangle} +
\Tr\left(\mathit{\Lambda}^\tau
T\right)\right|^2\leq
\left[\alpha a+\Tr\left(\mathit{\Lambda}^\tau
B\right)+2\Re\langle C^*,\mathit{\Gamma}^\tau
\rangle\right]\Tr\left(\mathit{\Lambda}^\tau U\right)
\end{equation}
holds for every $\alpha\in\bC$, matrices
$\mathit{\Gamma}=\left[\begin{array}{cccc}\gamma_1&\gamma_2&\ldots&\gamma_n\end{array}\right]$
and
$\mathit{\Lambda}=\left[\begin{array}{cccc}\lambda_{11}&\lambda_{12}&\ldots&\lambda_{1n}\\
\lambda_{21}&\lambda_{22}&\ldots&\lambda_{2n}\\\vdots&\vdots&
&\vdots\\\lambda_{n1}&\lambda_{n2}&\ldots&\lambda_{nn}\end{array}\right]$,
$\gamma_i\in\bC$, $\lambda_{ij}\in\bC$ for $i,j=1,2,\ldots,n$, such
that
\begin{enumerate}
\item $\alpha\geq 0$ and $\mathit{\Lambda}\geq 0$,
\item $\mathit{\Gamma}^*\mathit{\Gamma}\leq\alpha\mathit{\Lambda}$.
\end{enumerate}
The superscript $\tau$ denotes the transposition of matrices.
\end{lemma}
\begin{proof}
Obviously, the map $\varphi$ is positive if and only if
$\omega\circ\varphi$ is a positive functional on $M_2(\bC)$ for
every positive functional $\omega$ on $M_{n+1}(\bC)$.

Let $\omega$ be a linear functional on $M_{n+1}(\bC)$. Recall that
positivity of $\omega$ is equivalent to its complete positivity.
Hence, the Choi matrix
$\mathbf{H}_\omega=[\omega(F_{ij})]_{i,j=1}^{n+1}$ of $\omega$ is a
positive element of $M_{n+1}(\bC)$. Let us denote
$\alpha=\omega(F_{11})$, $\gamma_j=\omega(F_{1,j+1})$ for
$j=1,2,\ldots,n$, $\lambda_{ij}=\omega(F_{i+1,j+1})$ for
$i,j=1,2,\ldots,n$ and $\mathit{\Gamma}$ and $\mathit{\Lambda}$ are
defined as in the statement of the lemma. Then, we came to the
conclusion that $\omega$ is a positive functional if and only if the
matrix
\begin{equation}\label{matr3}
\left[\begin{array}{cc}\alpha&\mathit{\Gamma}\\\mathit{\Gamma}^*&\mathit{\Lambda}\end{array}\right]
\end{equation}
is a positive element of $M_{n+1}(\bC)$. This is equivalent to the
conditions (1) and (2) from the statement of the lemma.

Similarly, if $\omega'$ is a linear functional on $M_2(\bC)$, then
its positivity is equivalent to positivity of the matrix
$[\omega'(E_{ij})]_{i,j=1}^2$. Consequently, $\omega'$ is positive
if and only if $\omega'(E_{ii})\geq 0$ for $i=1,2$,
$\omega'(E_{21})=\ov{\omega'(E_{12})}$ and
\begin{equation}\label{om}
\left|\omega'(E_{12})\right|^2\leq
\omega'(E_{11})\omega'(E_{22})\end{equation} (cf. Corollary 8.4 in
\cite{St}).

Now, assume that $\alpha$, $\mathit{\Gamma}$ and $\mathit{\Lambda}$
are given and the conditions (1) and (2) are fulfilled.
In described above way it corresponds to some positive functional
$\omega$ on $M_{n+1}(\bC)$. Let $\omega'=\omega\circ\varphi$. Then
\begin{eqnarray}
\omega'(E_{11})&=& \alpha a + \sum_{i,j=1}^n\lambda_{ij}b_{ij}
+2\Re\left(\sum_{i=1}^n\gamma_ic_i\right)=\alpha
a+\Tr(\mathit{\Lambda}^\tau B)+2\Re\langle C^*,\mathit{\Gamma}^\tau
\rangle\\
\omega'(E_{22})&=&
\sum_{i,j=1}^n\lambda_{ij}u_{ij}=\Tr(\mathit{\Lambda}^\tau U)\\
\omega'(E_{12})&=&
\sum_{i=1}^n\gamma_iy_i+\sum_{i=1}^n\ov{\gamma_iz_i}+\sum_{i,j=1}^n\lambda_{ij}t_{ij}
=\langle Y^*,\mathit{\Gamma}^\tau\rangle+\ov{\langle Z^*,\mathit{\Gamma}^\tau
\rangle}+\Tr(\mathit{\Lambda}^\tau T).
\end{eqnarray}
It follows from the above equalities that (\ref{om}) is equivalent
to the inequality (\ref{nier1}). Hence, the statement of the lemma
follows from the remark contained in the first paragraph of the
proof.
\end{proof}

\begin{proposition}\label{Cmyz}
If the assumptions of Propositions \ref{Cmpos} are fulfilled, then
\begin{equation}\label{yz}
|Y|+|Z|\leq a^{1/2}U^{1/2}.
\end{equation}
\end{proposition}
\begin{proof}
Firstly, let us observe that the inequality (\ref{nier1}) can be
written in the form
\begin{eqnarray*}
\lefteqn{\left|\langle Y^*,\mathit{\Gamma}^\tau\rangle+\ov{\langle Z^*,\mathit{\Gamma}^\tau
\rangle}\right|^2+\left|\Tr\left(\mathit{\Lambda}^\tau
T\right)\right|^2+}\\
&+&2\Re\left[\left(\langle Y^*,\mathit{\Gamma}^\tau\rangle+\ov{\langle Z^*,\mathit{\Gamma}^\tau
\rangle}\right)\ov{\Tr\left(\mathit{\Lambda}^\tau
T\right)}\right]\leq
\left[\alpha a+\Tr\left(\mathit{\Lambda}^\tau
B\right)+2\Re\langle C^*,\mathit{\Gamma}^\tau
\rangle\right]\Tr\left(\mathit{\Lambda}^\tau
U\right)
\end{eqnarray*}
Putting $-\mathit{\Gamma}$ instead of $\mathit{\Gamma}$ we preserve
the positivity of the matrix (\ref{matr3}) and the above inequality
takes the form
\begin{eqnarray*}
\lefteqn{\left|\langle Y^*,\mathit{\Gamma}^\tau\rangle+\ov{\langle Z^*,\mathit{\Gamma}^\tau
\rangle}\right|^2+\left|\Tr\left(\mathit{\Lambda}^\tau
T\right)\right|^2+}\\
&-&2\Re\left[\left(\langle Y^*,\mathit{\Gamma}^\tau
\rangle+\ov{\langle Z^*,\mathit{\Gamma}^\tau
\rangle}\right)\ov{\Tr\left(\mathit{\Lambda}^\tau
T\right)}\right]\leq
\left[\alpha a+\Tr\left(\mathit{\Lambda}^\tau
B\right)-2\Re\langle C^*,\mathit{\Gamma}^\tau
\rangle\right]\Tr\left(\mathit{\Lambda}^\tau
U\right)
\end{eqnarray*}
If we add both above inequalities and divide the result by $2$, then
we get
\begin{equation*}
\left|\langle Y^*,\mathit{\Gamma}^\tau
\rangle+\ov{\langle Z^*,\mathit{\Gamma}^\tau
\rangle}\right|^2+\left|\Tr\left(\mathit{\Lambda}^\tau
T\right)\right|^2\leq\left[\alpha a+\Tr\left(\mathit{\Lambda}^\tau
B\right)\right]\Tr\left(\mathit{\Lambda}^\tau U\right)
\end{equation*}
and consequently
\begin{eqnarray}
\lefteqn{\left|\langle Y^*,\mathit{\Gamma}^\tau\rangle\right|^2+
\left|\langle Z^*,\mathit{\Gamma}^\tau\rangle\right|^2+
2\Re\langle Y^*,\mathit{\Gamma}^\tau\rangle\langle Z^*,\mathit{\Gamma}^\tau\rangle\leq}\label{nier2}\\
&\leq&\alpha a\Tr(\mathit{\Lambda}^\tau U)+\Tr(\mathit{\Lambda}^\tau B)\Tr(\mathit{\Lambda}^\tau U)
-|\Tr(\mathit{\Lambda}^\tau T)|^2\nonumber
\end{eqnarray}

Now, let $\eta$ be an arbitrary unit vector from $\bC^n$.
Then
\begin{eqnarray}
\lefteqn{\Vert (|Y|+|Z|)\eta\Vert\leq}\label{nieryz}\\
&\leq& \Vert |Y|\eta\Vert +\Vert |Z|\eta\Vert=
\Vert Y\Vert\,\Vert P_{\xi_Y}\eta\Vert + \Vert Z\Vert\,\Vert P_{\xi_Z}\eta\Vert\nonumber\\
&=&
\Vert Y\Vert\,|\langle\xi_Y,\eta\rangle|+\Vert Z\Vert\,|\langle\xi_Z,\eta\rangle|=
|\langle Y^*,\eta\rangle|+|\langle Z^*,\eta\rangle|\nonumber
\end{eqnarray}

Let $\varepsilon>0$. Put $\alpha=\varepsilon^{-1}$ and let $\mathit{\Gamma}=e^{i\theta}\eta^\tau$,
where $\theta$ is such that 
$$e^{2i\theta}\langle Y^*,\eta\rangle\langle Z^*,\eta\rangle=|\langle Y^*,\eta\rangle| |\langle Z^*,\eta\rangle|$$
hence
$$\langle Y^*,\mathit{\Gamma}^\tau\rangle\langle Z^*,\mathit{\Gamma}^\tau\rangle=
|\langle Y^*,\mathit{\Gamma}^\tau\rangle| |\langle Z^*,\mathit{\Gamma}^\tau\rangle|.$$
Put also $\mathit{\Lambda}=\varepsilon\mathit{\Gamma}^*\mathit{\Gamma}$.
Then from (\ref{nier2}) and (\ref{nieryz}) we have
$$
\Vert (|Y|+|Z|)\eta\Vert^2\leq
a\Tr((\mathit{\Gamma}^*\mathit{\Gamma})^\tau U)+\varepsilon^2
\left(\Tr(\mathit{\Lambda}^\tau B)\Tr(\mathit{\Lambda}^\tau U)
-|\Tr(\mathit{\Lambda}^\tau T)|^2\right)
$$
Since $\varepsilon$ is arbitrary then we have
$$
\Vert (|Y|+|Z|)\eta\Vert^2\leq
a\Tr((\mathit{\Gamma}^*\mathit{\Gamma})^\tau U)=
a\Tr(P_\eta U)=a\langle\eta,U\eta\rangle=\Vert a^{1/2}U^{1/2}\eta\Vert^2$$
and the proof is finished since $(|Y|+|Z|)^2\leq aU$ implies $|Y|+|Z|\leq a^{1/2}U^{1/2}$.
\end{proof}

\begin{lemma}\label{CP}
Let $\varphi:M_2(\bC)\to M_{n+1}(\bC)$ be a linear map with the Choi
matrix of the form (\ref{Cm3}). Then
\begin{enumerate}
\item
the map $\varphi$ is completely positive if and only if the
following conditions hold:
\begin{enumerate}
\item[(A1)] $Z=0$,
\item[(A2)] the matrix
$\left[\begin{array}{ccc}a&C&Y\\C^*&B&T\\Y^*&T^*&U\end{array}\right]$
is a positive element of $M_{2n+1}(\bC)$.
\end{enumerate}
In particular, the condition (A2) implies:
\begin{enumerate}
\item[(A3)] if $B$ is an invertible matrix, then $T^*B^{-1}T\leq U$,
\item[(A4)] $C^*C\leq aB$,
\item[(A5)] $Y^*Y\leq aU$,
\end{enumerate}
\item
the map $\varphi$ is completely copositive if and only if
\begin{enumerate}
\item[(B1)] $Y=0$,
\item[(B2)] the matrix
$\left[\begin{array}{ccc}a&C&Z\\C^*&B&T^*\\Z^*&T&U\end{array}\right]$
is a positive element of $M_{2n+1}(\bC)$.
\end{enumerate}
In particular, the condition (B2) implies:
\begin{enumerate}
\item[(B3)] if $B$ is invertible, then $TB^{-1}T^*\leq U$,
\item[(B4)] $C^*C\leq aB$,
\item[(B5)] $Z^*Z\leq aU$,
\end{enumerate}
\end{enumerate}
\end{lemma}
\begin{proof}
It is rather obvious that the conditions (A1) and (A2) imply
positivity of the matrix (\ref{Cm3}), and consequently the complete
positivity of the map $\varphi$. On the other hand it is easy to see
that (A2) is a necessary condition for positivity of the matrix
(\ref{Cm3}). In order to finish the proof of the first part of point
(1) one should show that positivity of (\ref{Cm3}) implies $Z=0$.

Let $L_1$ be a linear subspace generated by the vector $f_1$ and let
$L_2$ be a subspace spanned by $f_2,f_3,\ldots,f_{n+1}$, so $\bC^{n+1}=L_1\oplus
L_2$. Any vector $v\in\bC^{n+1}$ can be uniquely decomposed onto the sum
$v=v^{(1)}+v^{(2)}$, where $v^{(i)}\in L_i$, $i=1,2$. Blocks of the
matrix (\ref{Cm3}) are interpreted as operators. Namely:
$B,T,U:L_2\to L_2$, $C,Y,Z:L_2\to L_1$, and $a:L_1\to L_1$. Recall
(cf. \cite{T}), that the positivity of the matrix (\ref{Cm3}) is
equivalent to the following inequality
\begin{equation*}
\left\langle
v_1,\left[\begin{array}{cc}a&C\\C^*&B\end{array}\right]v_1\right\rangle+
\left\langle
v_1,\left[\begin{array}{cc}0&Y\\Z^*&T\end{array}\right]v_2\right\rangle
+ \left\langle v_2,\left[\begin{array}{cc}0&Z
\\ Y^*&T^*\end{array}\right]v_1\right\rangle + \left\langle
v_2,\left[\begin{array}{cc}0&0
\\ 0&U\end{array}\right]v_2\right\rangle\geq 0
\end{equation*}
for any $v_1,v_2\in\bC^{n+1}$. This is equivalent to
\begin{eqnarray}\label{CP1}
\lefteqn{\left\langle v_1^{(1)},av_1^{(1)}\right\rangle+
\left\langle
v_1^{(2)},Bv_1^{(2)}\right\rangle+ \left\langle v_2^{(2)},Uv_2^{(2)}\right\rangle+}\\
&+&2\Re\left\langle v_1^{(1)},Cv_1^{(2)}\right\rangle+
2\Re\left\langle v_1^{(1)},Yv_2^{(2)}\right\rangle+ 2\Re\left\langle
v_2^{(1)},Zv_1^{(2)}\right\rangle+ 2\Re\left\langle
v_1^{(2)},Tv_2^{(2)}\right\rangle\geq 0\nonumber
\end{eqnarray}
where $v_j=v_j^{(1)}+v_j^{(2)}$ for $j=1,2$, and
$v_1^{(1)},v_2^{(1)}\in L_1$ and $v_1^{(2)},v_2^{(2)}\in L_2$.
Assume that $v_1^{(1)}=0$, $v_2^{(2)}=0$, $v_1^{(2)}$ is an
arbitrary element of $L_2$, and $v_2^{(1)}=-rZv_1^{(2)}$ for some
$r>0$. Then (\ref{CP1}) reduces to
\begin{equation}
\left\langle v_1^{(2)},Bv_1^{(2)}\right\rangle -r\left\Vert
Zv_1^{(2)}\right\Vert^2\geq 0.
\end{equation}
This inequality holds for any $v_1^{(2)}\in L_2$ and $r>0$. It is
possible only for $Z=0$.

To show the second part of the point (1) one needs to notice that
positivity of the matrix
$\left[\begin{array}{ccc}a&C&Y\\C^*&B&T\\Y^*&T^*&U\end{array}\right]$
implies the positivity of the matrices
$\left[\begin{array}{cc}B&T\\T^*&U\end{array}\right]$,
$\left[\begin{array}{cc}a&C\\C^*&B\end{array}\right]$ and
$\left[\begin{array}{cc}a&Y\\Y^*&U\end{array}\right]$.

Recall that the map is completely copositive if and only if the
partial transposition of the matrix (\ref{Cm3}) is positive. The
partial transposition (cf. \cite{MM2}) of (\ref{Cm3}) is equal to
$$\left[\begin{array}{cc|cc}a&C&0&Z\\C^*&B&Y^*&T^*\\\hline
0&Y&0&0\\Z^*&T&0&U\end{array}\right].$$ So, to prove the point (2)
we can use the arguments as in proof of the point (1).
\end{proof}

Now, assume that $\varphi:M_2(\bC)\to M_3(\bC)$ is a unital positive
map, and $\varphi\in F_{e_2,f_1}$. Hence its Choi matrix has the
form
\begin{equation}\label{Cmunit}
\left[\begin{array}{cc|cc} 1&0&0&Y\\0&B&Z^*&T\\\hline
0&Z&0&0\\Y^*&T^*&0&U\end{array}\right],\end{equation} where $B$ and
$U$ are positive matrices such that $B+U=1$ and conditions listed in
Propositions \ref{Cmpos} and \ref{Cmyz} are satisfied. From the
theorem of Woronowicz (cf. \cite{Wor}) it follows that there are maps
$\varphi_1,\varphi_2:M_2(\bC)\to M_3(\bC)$ such that
$\varphi=\varphi_1+\varphi_2$, and $\varphi_1$ is a completely
positive map while $\varphi_2$ is a completely copositive one. From
Definition \ref{face} we conclude that both $\varphi_1$ and
$\varphi_2$ are contained in the face $F_{e_2,f_1}$. So, from
Proposition \ref{CP} it follows that their Choi matrices
$\mathbf{H}_1$ and $\mathbf{H}_2$ are of the form
\begin{equation}
\mathbf{H}_1=\left[\begin{array}{cc|cc}a_1&C&0&Y\\C^*&B_1&0&T_1\\\hline
0&0&0&0\\Y^*&T_1^*&0&U_1\end{array}\right],\qquad
\mathbf{H}_2=\left[\begin{array}{cc|cc}a_2&-C&0&0\\-C^*&B_2&Z^*&T_2\\\hline
0&Z&0&0\\0&T_2^*&0&U_2\end{array}\right],
\end{equation}
where $a_i,B_i,U_i\geq 0$ for $i=1,2$, and the following equalities
hold: $a_1+a_2=a$, $T_1+T_2=T$, $B_1+B_2=B$ and $U_1+U_2=U$. In
\cite{MM2} we proved that if $\varphi:M_2(\bC)\to M_2(\bC)$ is from
a large class of extremal positive unital maps, then the maps
$\varphi_1$ and $\varphi_2$ are uniquely determined (cf. Theorem 2.7
in \cite{MM2}). Motivated by this type of decomposition and the results given in this section 
(we `quantized' the relations (1)-(3) given at the end of Section 1) 
we wish to formulate the
following conjecture:

\textit{Assume that $\varphi:M_2(\bC)\to M_3(\bC)$ is a positive
unital map such that $U\neq 0$, $Y\neq 0$, $Z\neq 0$ and
$|Y|+|Z|=U^{1/2}$. Then the
decomposition $\varphi=\varphi_1+\varphi_2$ onto completely positive
and completely copositive parts is uniquely determined.}

\end{document}